\def \NN {\mathbb N}
\def \CC {\mathbb C}

\def \RR {\mathbb R}

\def \epsilon{\varepsilon}

\def \d {\text{d}}

\def \fine {{\hfill \qedsymbol}}

\def \si {\sigma}

\documentclass[12pt,reqno]{amsart}

\usepackage{amsmath,amssymb}
\DeclareMathOperator{\res}{res}

\numberwithin{equation}{section}

\begin{document}

\title[]{Explicit formulae for averages \\ of Goldbach representations}

\author[]{J.Br\"udern, J.Kaczorowski \lowercase{and} A.Perelli}

\date{}
\maketitle

\bigskip
{\bf Abstract.} We prove an explicit formula, analogous to the classical explicit formula for $\psi(x)$, for the Ces\`aro-Riesz mean of any order $k>0$ of the number of representations of $n$ as a sum of two primes. Our approach is based on a double Mellin transform and the analytic continuation of certain functions arising therein.



\vskip1cm
\section{Introduction}

\smallskip
Riemann's explicit formula reveals an intimate relation between the distribution of primes and the zeros of the Riemann zeta function. A similar formula for the number of prime points in a triangle such as
\begin{equation}
\label{1-1}
\sharp\{(p,p'): p+p'\leq N, \ \text{$p$ and $p'$ prime}\}
\end{equation}
may be expected, and Fujii \cite{Fuj/1991} was first to step in this direction. In recent years, there has been a flurry of papers on the topic,  also analyzing weighted versions of the count \eqref{1-1}; see Bhowmik \& Schlage-Puchta \cite{Bo-SP/2010}, Languasco \& Zaccagnini \cite{La-Za/2012},\cite{La-Za/2015}, Goldston \& Yang \cite{Go-Ya/arxiv} and Languasco's survey \cite{Lan/2016}. So far, all authors have interpreted the quantity in \eqref{1-1} as the mean of the number of representations of a given number $n\leq N$ as the sum of two primes, and replaced the latter by the analytically simpler expression
\[
R(n) = \sum_{m+m'=n}\Lambda(m)\Lambda(m'),
\]
with $\Lambda$ the von Mangoldt function. One is then led to study the sum
\begin{equation}
\label{1-2}
G_0(N) = \sideset{}{'} \sum_{n\leq N} R(n),
\end{equation}
where the notation indicates that $R(N)/2$ is to be subtracted from the sum in \eqref{1-2} should $N$ be a natural number.

\smallskip
A purely formal application of the method to be described herein suggests that, perhaps, one has 
\begin{equation}
\label{1-3}
G_0(N) = \lim_{T\to\infty} \sum_{\substack{|\sigma|\leq T,|t|\leq T \\ |u|\leq T,|v|\leq T}}  \underset{s}{\res}\, \underset{w}{\res}\, \frac{\zeta'}{\zeta}(s) \frac{\zeta'}{\zeta}(w)    N^{s+w} \frac{\Gamma(s) \Gamma(w)}{\Gamma(w+s+1)}  .
\end{equation}
Here and later, we write $s=\sigma+it$, $w=u+iv$, with real numbers $\sigma,t,u,v$. If true, then \eqref{1-3}  would correspond to Riemann's formula, in the form
\[
 \sideset{}{'} \sum_{n\leq N} \Lambda(n) = -\lim_{T\to\infty} \sum_{|\sigma|\leq T,|t|\leq T}\underset{s}{\res} \, \frac{\zeta'}{\zeta}(s) \frac{N^s}{s}
\]
(see Chapter 17 of Davenport \cite{Dav/1980}), in a two-dimensional context.

\smallskip
Our purpose is to present an approach that yields  weighted versions of \eqref{1-3} superior to all previous work. {In much recent work, } $G_0(N)$ was examined by means of the circle method, but this procedure apparently puts limitations to the insights obtainable on the problem. Our method is radically different; indeed, we consider the counting problem in its natural framework, as a two-dimensional lattice point count with constraints on the coordinates. The range of summation $m+m'\leq N$ constitutes a triangle, and the Mellin transform of its indicator function is the analytic function in $s$ and $w$ that occurs on the right hand side of \eqref{1-3}.  In this way, our hypothetical formula \eqref{1-3} arises canonically, and any attempt to establish it may well lead to approximations to the conjecture in full.

\smallskip
 In a weighted setting,  the ideas described above take us well beyond the current state of affairs. We follow Languasco \& Zaccagnini \cite{La-Za/2015} and consider the Ces\`aro-Riesz mean of $R(n)$ of order $k>0$, defined for $N\geq 2$ by
\begin{equation}
\label{1-4}
G_k(N) = \frac{1}{\Gamma(k+1)} \sum_{n< N}R(n) \left(1-\frac{n}{N}\right)^k.
\end{equation}
They show {(\cite[Theorem 1]{La-Za/2015}, corrected in \cite{Lan/2016}) } that whenever $k>1$ one has
\begin{equation}
\label{1-5}
G_k(N) = \frac{N^2}{\Gamma(k+3)} -2A_k(N) + B_k(N) + O(N),
\end{equation}
in which
\begin{equation}
\label{1-6}
A_k(N) = \sum_{\rho} \frac{\Gamma(\rho)}{\Gamma(\rho+k+2)} N^{\rho+1} \quad \text{and} \quad  B_k(N) = \sum_{\rho}\sum_{\rho'} \frac{\Gamma(\rho)\Gamma(\rho')}{\Gamma(\rho+\rho'+k+1)} N^{\rho+\rho'}
\end{equation}
with the sums running  over the non-trivial zeros of $\zeta(s)$, each zero listed according to its multiplicity. We apply this notation throughout the paper.
 Note that the sums in \eqref{1-6} are absolutely convergent for $k>1/2$. Furthermore, observe that  \eqref{1-5} is compatible with a weighted version of \eqref{1-3} in which only the singularities with $\sigma>0$ and $u>0$ are made explicit. 
We are able to handle the wider range $k>0$ while
taking into account all the singularities, hence  providing a full explicit formula for $G_k(N)$. This improves on 
Languasco \& Zaccagnini \cite{La-Za/2015} in several directions,
and a  recent result by Goldston \& Yang \cite{Go-Ya/arxiv} for $k=1$, obtained under the Riemann Hypothesis, is also superseded unconditionally.

\smallskip
To initiate our treatment, 
let $N\geq 4$ be an integer, let $k>0$ be real, and then  observe that
\[
1- \frac{m+n}{N} = \Big(1- \frac{n}{N-m}\Big)  \Big(1- \frac{m}{N}\Big)
\]
to recast \eqref{1-4} as
\begin{equation}
\label{3-1}
G_k(N) = \frac{1}{\Gamma(k+1)}  \sum_{m< N}  \Lambda(m) \left(1-\frac{m}{N}\right)^k  \sum_{n< N-m}  \Lambda(n) \left(1-\frac{n}{N-m}\right)^k.
\end{equation}
We now recall that whenever $x$ and $c$ are positive numbers, and $z$ is a complex number with  
with $\Re z>0$, then
\begin{equation}\label{mellin}
\frac{1}{2\pi i} \int_{(c)} \frac{\Gamma(s) x^{-s}}{\Gamma(s+z+1)}\, \d s =
\begin{cases}
(1-x)^z/\Gamma(z+1) & \text{if} \ 0<x<1, \\
0 & \text{if} \ x\geq 1
\end{cases}
\end{equation}
(for example, this is formula 5.35 of Section 2.5 on p.195 of Oberhettinger \cite{Obe/1974}, with the choice $\alpha=0$ and $\beta = z+1$). We apply \eqref{mellin} twice, first with $z=k$ to the inner sum in \eqref{3-1} and then with $z=s+k$ to the outer sum. 
This leads us to
\begin{equation}
\label{3-2}
G_k(N) = \frac{1}{(2\pi i)^2} \int_{(2)} \int_{(2)} \frac{\zeta'}{\zeta}(w)  \frac{\zeta'}{\zeta}(s) \frac{\Gamma(w)\Gamma(s)}{\Gamma(s+w+k+1)} N^{w+s} \,\d s\, \d w.
\end{equation} 

\smallskip
To proceed further, we shift the inner integration from the line $\Re s = 2$ to $\Re s = {-} 1/2$. We shall then encounter  the functions 
\begin{equation}
\label{1-7}
T_N(w) = \frac{{-}1}{2\pi i} \int_{(-1/2)} \frac{\zeta'}{\zeta}(s) \frac{\Gamma(s)}{\Gamma(s+w+1)} N^s\, \d s
\end{equation}
and
\begin{equation}
\label{1-8}
Z_N(w) = \sum_{\rho} \frac{\Gamma(\rho)}{\Gamma(\rho+w+1)} N^\rho.
\end{equation}
By  Stirling's formula, the integral  in \eqref{1-7} and the sum in \eqref{1-8} are absolutely and compactly  convergent  in $u>0$, hence $T_N(w)$ and $Z_N(w)$ are holomorphic in this half-plane. It is fundamental to this paper to realize that
the functions $T_N(w)$ and $Z_N(w)$ both extend to entire functions that do not grow too fast. We establish this in Section 2. Equipped with this, it will be possible to shift the integration over $\Re w = 2$ to $\Re w = - M$, say, and to analyse the limit
for $M\to \infty$. We then arrive at the desired explicit formula for $G_k(N)$. The final result features  the sums
\begin{equation}
\label{1-12}
\Sigma_\Gamma(N,k) = - \sum_{\nu=1}^\infty \underset{w=-\nu}{\res}\, \frac{\zeta'}{\zeta}(w) \Gamma(w) \frac{N^w}{\Gamma(w+k+1)} ,
\end{equation}
\begin{equation}
\label{1-13}
\Sigma_Z(N,k) = -\sum_{\nu=1}^\infty \underset{w=-\nu}{\res}\, \frac{\zeta'}{\zeta}(w) \Gamma(w) Z_N(w+k) N^w ,
\end{equation}
\begin{equation}
\label{1-14}
\Sigma_T(N,k) = -\sum_{\nu=1}^\infty \underset{w=-\nu}{\res}\,\frac{\zeta'}{\zeta}(w) \Gamma(w) T_N(w+k) N^w .
\end{equation}

\medskip
{\bf Theorem.} {\sl Let $N\geq4$ be an integer and $k>0$. Then}
\[
\begin{split}
G_k(N) &= \frac{N^2}{\Gamma(k+3)} - 2NZ_N(k+1) + \sum_\rho \Gamma(\rho) Z_N(\rho+k) N^\rho -2 \frac{\zeta'}{\zeta}(0) \frac{N}{\Gamma(k+2)} \\
&\hskip-.7cm + 2\frac{\zeta'}{\zeta}(0) Z_N(k) + NT_N(k+1) + \frac{\zeta'}{\zeta}(0)^2 \frac{1}{\Gamma(k+1)} - \sum_\rho \Gamma(\rho) T_N(\rho+k) N^\rho - \frac{\zeta'}{\zeta}(0) T_N(k) \\
&  + N\Sigma_\Gamma(N,k+1) - \Sigma_Z(N,k) - \frac{\zeta'}{\zeta}(0) \Sigma_\Gamma(N,k) + \Sigma_T(N,k),
\end{split}
\]
{\sl where the sums \eqref{1-12}, \eqref{1-13}, \eqref{1-14} and the sums  over the non-trivial zeros of $\zeta(s)$ are absolutely convergent.} 

\medskip
Some remarks are in order. First note that 
\[
NZ_N(k+1) = A_k(N) \ \text{for} \ k>0 \quad \text{and}  \quad \sum_\rho \Gamma(\rho) Z_N(\rho+k) N^\rho = B_k(N) \ \text{for} \ k>1/2,
\]
so that  the first three terms in our explicit formula correspond to the explicit terms in \eqref{1-5}. The other terms are listed, roughly, in descending order according to their expected order of magnitude in  $N$. Moreover, the first nine terms in the explicit formula correspond to the contribution of the poles of $-\Gamma(s)\zeta'(s)/\zeta(s)$  in $\sigma\geq0$, see \eqref{3-6}--\eqref{3-8} below. The explicit form of the sums in \eqref{1-12}, \eqref{1-13} and \eqref{1-14}, corresponding to the contribution of the remaining poles, can be found in \eqref{3-16}, \eqref{3-17} and \eqref{3-18} below. We also remark that the series expansions for $\Sigma_\Gamma(N,k)$ and $\Sigma_T(N,k)$ in \eqref{3-16} and \eqref{3-18} are actually asymptotic expansions. Indeed, \eqref{3-19} below shows that cutting these  series at $\nu=M\geq 2$ produces an error of size, roughly, $O(N^{-(M+1)})$. This does not hold for the sum $\Sigma_Z(N,k)$, in which case we can only prove that the tail is of order $O(N^{-k+\epsilon})$ for every $\epsilon>0$; see again \eqref{3-19}.

\smallskip
Thus far, we have concentrated on a full explicit formula, but it is worth pointing out that truncated versions are also available, not necessarily obtained by cutting the formula in our theorem. It may be instructive to compare the potential of our approach with that used by Languasco and Zaccagnini \cite{La-Za/2015}. As we have pointed out already, in order to identify the third term in our formula with the sum
$B_k(N)$ we require that $k>1/2$. Given this condition, it is then immediate from \eqref{1-7} and \eqref{1-8}
that \eqref{1-5} remains valid in this range for $k$, but of course one can replace the $O(N)$ term by a more explicit expression. For example, if one uses the zero-free region for $\zeta(s)$ within \eqref{1-8}, one readily confirms that
\[
\begin{split}
G_k(N)& = \frac{N^2}{\Gamma(k+3)} - 2A_k(N) + B_k(N) 
 -2 \frac{\zeta'}{\zeta}(0) \frac{N}{\Gamma(k+2)}+ O(N \exp(-c \sqrt{\log N}),
\end{split}
\]
with some $c>0$. However, a sharper truncated formula can be obtained by shifting the integrations in \eqref{3-2} to the left as much as is allowed by the convergence of the involved quantities. For example, shifting the integrations to (roughly) $-3/2$, we infer that in the range $k>1/2$ one has
\[
\begin{split}
G_k(N) = \frac{N^2}{\Gamma(k+3)} - 2A_k(N) + B_k(N)  
-2 \frac{\zeta'}{\zeta}(0)\frac{N}{\Gamma(k+2)}  +2\frac{\zeta'}{\zeta}(0) Z_N(k) + C + o(1)
\end{split}
\]
with
\[
C =  \left( \frac{\zeta'}{\zeta}(0)^2+ 2\frac{\zeta'}{\zeta}(-1)\right) \frac{1}{\Gamma(k+1)},
\]
which, up to an error of size $o(1)$, expresses 
$G_k(N)$ in terms  of $A_k$, $B_k$ and $Z_N(k)$, all directly defined by sums of the zeros of the Riemann zeta function.

\bigskip
{\bf Acknowledgements.} This research was partially supported by a grant from Deutsche Forschungsgemeinschaft, the grant PRIN-2015 {\sl ``Number Theory and Arithmetic Geometry''} and by a grant from the 
National Science Centre, Poland; AP is member of the GNAMPA group of INdAM.

\bigskip
\section{The auxiliary functions}

\smallskip
{\bf 2.1. The analytic continuation.} In this section, we are concerned with the auxiliary functions 
$T_N(w)$ and $Z_N(w)$. Most of our efforts go into the proof of the following facts.

\medskip
{\bf Proposition 1.} {\sl Let $N\geq4$, and let $\gamma$ denote Euler's constant. Then $T_N(w)$ extends to an entire function and satisfies}
\begin{equation}
\label{1-9}
\begin{split}
T_N(w)  = &\frac{1}{\Gamma(w+1)} \Big\{ -\sum_{n=1}^\infty \frac{\Lambda(n)}{n} \Big(\big(1-\frac{1}{nN}\big)^w-1\Big) + \int_0^1 \Big(\big(1-\frac{\xi}{N}\big)^w-1\Big)\, \frac{\d \xi}{\xi} \\
& - \Big(\big(1-\frac{1}{N}\big)^w-1\Big) \log(2\pi e^\gamma) - { \frac{1}{2} \big(1-\frac{1}{N}\big)^w } \log \big(1-\frac{1}{N^2}\big) \\
& - \int_1^N \Big(\big(1-\frac{1}{\xi}\big)^w - \big(1-\frac{1}{N}\big)^w \Big) \frac{\xi\,\d \xi}{N^2-\xi^2} \\
& + N\int_N^\infty \Big(\big(1-\frac{1}{N}\big)^w - \big(1-\frac{1}{\xi}\big)^w \Big) \frac{\d \xi}{\xi^2-N^2} \Big\}.
\end{split}
\end{equation}
{\sl Moreover, there is a real number $K$ such that  for any $\delta$ with $0<\delta <1$ and $|w+m|>\delta$ for all integers $m\geq1$ we have}
\begin{equation}
\label{1-10}
T_N(w) \le K  \frac{2^{|u|}\log(|w|+2)}{\delta |\Gamma(w+1)|}.
\end{equation}

\medskip
There is a similar yet less explicit result for the function $Z_N(w)$.

\medskip
{\bf Proposition 2.} {\sl Let $N\geq 4$. Then $Z_N(w)$ extends to an entire function. Moreover,  there is a real number $K$ such that  for any $\delta$ with $0<\delta<1$ and $|w+m|>\delta$ for all integers $m\geq1$ we have
\begin{equation}
\label{1-11}
Z_N(w) \le \frac{ K}{\delta |\Gamma(w+1)|} \times
\begin{cases}
\big( N^{|u|+1} +  2^{|u|}\log(|w|+2)\big) & \ \text{if $u\in\RR$}, \\
\big( N^{|u|}\log N + 2^{|u|}\log |w|\big) & \ \text{if $u\leq -3/2$}.
\end{cases}
\end{equation}
}

\medskip
{\bf 2.2. The function $T_N(w)$. First steps.} We begin with the proof of Proposition 1. 
From the functional equation of $\zeta(s)$ in the form
\[
\zeta(1-s) = 2 (2\pi)^{-s} \cos\big(\frac{\pi s}{2}\big) \Gamma(s) \zeta(s)
\]
we obtain
\begin{equation}
\label{2-1}
\frac{\zeta'}{\zeta}(s) = G(s) - \frac{\zeta'}{\zeta}(1-s),
\end{equation}
where
\begin{equation}
\label{2-2}
G(s) = \log(2\pi) - \frac{g'}{g}(s) \quad \text{and} \quad g(s) = \Gamma(s) \cos\frac{\pi s}{2}.
\end{equation}
Accordingly, for $N\geq 1$ and $u>0$ we have
\begin{equation}
\label{2-3}
\begin{split}
T_N(w) & =  -\frac{1}{2\pi i} \int_{(-1/2)} G(s) \frac{\Gamma(s)}{\Gamma(s+w+1)} N^s\, \d s \\
& \hskip.5cm + \frac{1}{2\pi i} \int_{(-1/2)} \frac{\zeta'}{\zeta}(1-s) \frac{\Gamma(s)}{\Gamma(s+w+1)} N^s \d s \\
&= -T_N^{(1)}(w) + T_N^{(2)}(w),
\end{split}
\end{equation}
say.
Since the integration is on the line $\si=-1/2$, we may expand $\frac{\zeta'}{\zeta}(1-s)$ and switch summation and integration, thus getting that
\begin{equation}
\label{2-4}
T_N^{(2)}(w) = -\sum_{n=1}^\infty \frac{\Lambda(n)}{n} \int_{(-1/2)} \frac{\Gamma(s)}{\Gamma(s+w+1)} (nN)^s \,\d s.
\end{equation}
We shift  the line of integration in \eqref{2-4} to $\si=c>0$ and compute  the residue at $s=0$. Recalling that $\Lambda(1)=0$, we infer from \eqref{mellin}  with $z=w$ and $x=(nN)^{-1}$ that
\[
T_N^{(2)}(w) =  \frac{{-}1}{\Gamma(w+1)}  \sum_{ n=1}^\infty  \frac{\Lambda(n)}{n} \Big(\big(1-\frac{1}{nN}\big)^w -1\Big).
\]
{Yet, } for $N\geq 1$, one has
\[
\big(1-\frac{1}{nN}\big)^w -1 \ll_{w} \frac{1}{n},
\]
{w}hence 
\begin{equation}
\label{2-5}
T_N^{(2)}(w) = \frac{{-}1}{\Gamma(w+1)} \sum_{n=1}^\infty \frac{\Lambda(n)}{n} \Big(\big(1-\frac{1}{nN}\big)^w -1\Big),
\end{equation}
the series being absolutely and compactly  convergent on $\CC$. In particular, for all $N\ge 1$, the function    $T_N^{(2)}(w)$ is entire. 

\medskip
{\bf 2.3. Computing $T_N^{(1)}(w)$.}  Recalling the definition of $g(s)$ in \eqref{2-2} and the Hadamard products of $1/\Gamma(s)$ and $\cos(\pi s/2)$ (Chapter 10 of Davenport \cite{Dav/1980} and Chapter 1 of Remmert \cite{Rem/1998}, respectively) we obtain
\[
g(s) = \frac{1}{s} e^{-\gamma s} \prod_{\substack{2|\nu \\ \nu\geq1}} \left(\frac{\nu}{s+\nu}\right) e^{s/\nu} \prod_{\substack{2\nmid \nu \\ \nu\geq1}} \left(1-\frac{s}{\nu}\right) e^{s/\nu},
\]
the products being absolutely and uniformly convergent on any compact subset of $\CC$ not containing a zero or pole of $g(s)$.  
Therefore, away from such points we have
\begin{equation}
\label{2-6}
\begin{split}
\frac{g'}{g}(s) &= - \frac{1}{s} - \gamma - \sum_{\substack{2|\nu \\ \nu\geq1}} \left(\frac{1}{s+\nu}- \frac{1}{\nu} \right) +  \sum_{\substack{2\nmid \nu \\ \nu\geq1}} \left(\frac{1}{s-\nu}+ \frac{1}{\nu} \right) \\
&= - \frac{1}{s} - \gamma - \Sigma_1(s) + \Sigma_2(s),
\end{split}
\end{equation}
say. Recalling \eqref{2-2} and inserting \eqref{2-6} into the expression for $T_N^{(1)}(w)$ in \eqref{2-3}, for $N\geq 1$ we get
\begin{equation}
\label{2-7}
\begin{split}
T_N^{(1)}(w) &= \frac{1}{2\pi i} \int_{(-1/2)} \frac{\Gamma(s)}{\Gamma(s+w+1)}\, \frac{N^s}{s} \,\d s 
 + \big(\log(2\pi) + \gamma \big) \frac{1}{2\pi i} \int_{(-1/2)} \frac{\Gamma(s)}{\Gamma(s+w+1)} N^s\, \d s \\
&\hskip.6cm + \frac{1}{2\pi i} \int_{(-1/2)} \Sigma_1(s) \frac{\Gamma(s)}{\Gamma(s+w+1)} N^s \,\d s - \frac{1}{2\pi i} \int_{(-1/2)} \Sigma_2(s) \frac{\Gamma(s)}{\Gamma(s+w+1)} N^s\,\d s \\
&= T_N^{(1,1)}(w) + \log\big(2\pi e^\gamma\big) T_N^{(1,2)}(w) + T_N^{(1,3)}(w) - T_N^{(1,4)}(w),
\end{split}
\end{equation}
say. Next we compute separately the functions $T_N^{(1,j)}(w)$, $j=1,\dots,4$, in \eqref{2-7}.

\medskip
{\bf 2.4. Computing $T_N^{(1,1)}(w)$ and $T_N^{(1,2)}(w)$.} We differentiate  $T_N^{(1,1)}(w)$ with respect to $N$ and then 
{apply \eqref{mellin} } with $z=w$ and $x=1/N$. For $N\geq 1$ and $u>0$ this yields 
\begin{equation}
\label{2-8}
\frac{\partial}{\partial N} T_N^{(1,1)}(w) = \frac{1}{2\pi iN} \int_{(-1/2)} \frac{\Gamma(s)}{\Gamma(s+w+1)} N^s \,\d s = \frac{1}{N\Gamma(w+1)} \Big(\big(1-\frac{1}{N}\big)^w-1\Big).
\end{equation}
We remark that, here and later, the expression $(1-\theta)^w$ equals 0 if $\theta=1$ and $u>0$. Moreover, from the expression of $T_N^{(1,1)}(w)$ in \eqref{2-7} and the first equality in \eqref{2-8} we see that 
\[
T_N^{(1,1)}(w) \ll N^{-1/2} \quad \text{and} \quad \frac{\partial}{\partial N} T_N^{(1,1)}(w) \ll N^{-3/2} \quad \text{as} \ N\to\infty,
\] 
hence from the second equality in \eqref{2-8} we get
\begin{equation}
\label{2-9}
T_N^{(1,1)}(w)  = -\int_N^\infty \frac{\partial}{\partial x} T_x^{(1,1)}(w)\, \d x = -\frac{1}{\Gamma(w+1)} \int_N^\infty  \Big(\big(1-\frac{1}{x}\big)^w-1\Big) \,\frac{\d x}{x}.
\end{equation}
Now we substitute
\[
x= \frac{N}{\xi}, \quad 0<\xi<1, \quad \frac{\d x}{x} = -\frac{\d \xi}{\xi}
\]
to see from \eqref{2-9} that
\begin{equation}
\label{2-10}
T_N^{(1,1)}(w)  = -\frac{1}{\Gamma(w+1)} \int_0^1  \Big(\big(1-\frac{\xi}{N}\big)^w-1\Big)\, \frac{\d \xi}{\xi}.
\end{equation}
But as $\xi\to 0^+$ we have $\big(1-\frac{\xi}{N}\big)^w-1 \ll \xi$, and this estimate holds uniformly as long as $w$ ranges over a fixed compact subset  of $\CC$. Hence \eqref{2-10} shows that $T_N^{(1,1)}(w)$ is an entire function for every $N\geq 2$.

\smallskip
We now turn to $T_N^{(1,2)}(w)$. Here, \eqref{mellin} immediately yields 
\begin{equation}
\label{2-11}
T_N^{(1,2)}(w)  = \frac{1}{\Gamma(w+1)}  \Big(\big(1-\frac{1}{N}\big)^w-1\Big),
\end{equation}
hence $T_N^{(1,2)}(w)$ is an entire function as well, for every $N\geq 2$.

\medskip 
{\bf 2.5. Computing $T_N^{(1,3)}(w)$.} Suppose that $N\geq 1$ and $u>0$. Then, recalling \eqref{2-6} and \eqref{2-7},  we infer from \eqref{mellin} that
\begin{equation}
\label{2-12}
\begin{split}
T_N^{(1,3)}(w)  &=  \sum_{\substack{\nu=2 \\ 2\mid \nu}}^\infty \Big\{ \frac{1}{2\pi i} \int_{(-1/2)} \frac{N^s}{s+\nu}\, \frac{\Gamma(s)}{\Gamma(s+w+1)}\, \d s 
 - \frac{1}{\nu \Gamma(w+1)}  \Big(\big(1-\frac{1}{N}\big)^w-1\Big) \Big\} \\
&=  \sum_{\substack{\nu=2 \\ 2\mid \nu}}^\infty \left\{ H_N(w,\nu)   - \frac{1}{\nu \Gamma(w+1)}  \Big(\big(1-\frac{1}{N}\big)^w-1\Big) \right\},
\end{split}
\end{equation}
say. We differentiate the identity
\begin{equation}
\label{2-13}
N^\nu H_N(w,\nu) = \frac{1}{2\pi i} \int_{(-1/2)} \frac{N^{s+\nu}}{s+\nu}\, \frac{\Gamma(s)}{\Gamma(s+w+1)} \d s,
\end{equation}
with respect to $N$ and then applying \eqref{mellin} to confirm  that $H_N(w,\nu)$ satisfies the differential equation
\begin{equation}
\label{2-14}
\frac{\partial}{\partial N} H_N(w,\nu) + \frac{\nu}{N} H_N(w,\nu) = \frac{1}{N\Gamma(w+1)} \Big(\big(1-\frac{1}{N}\big)^w-1\Big). 
\end{equation}
We solve \eqref{2-14} by searching for a function $c_N(w,\nu)$ such that
\begin{equation}
\label{2-15}
H_N(w,\nu) = c_N(w,\nu) N^{-\nu},
\end{equation}
which in view of \eqref{2-14} satisfies
\begin{equation}
\label{2-16}
\frac{\partial}{\partial N} c_N(w,\nu) = \frac{N^{\nu-1}}{\Gamma(w+1)} \Big(\big(1-\frac{1}{N}\big)^w-1\Big).
\end{equation}
Next we take $N=1$ in \eqref{2-13}. A computation based on Stirling's formula shows that for each $\nu\geq1$ one has
\[
\lim_{x\to+\infty} \frac{1}{2\pi i} \int_{(x)} \frac{1}{s+\nu}\, \frac{\Gamma(s)}{\Gamma(s+w+1)}\, \d s =0,
\]
{uniformly for $w$ in any fixed  compact part of  $u>0$}; we refer to Section 3 for a more explicit presentation of similar computations. Hence, from \eqref{2-13} and \eqref{2-15} with $N=1$, shifting the line of integration to $+\infty$ and computing the residue at $s=0$ we obtain the boundary condition
\[
c_1(w,\nu) = -\frac{1}{\nu \Gamma(w+1)}.
\]
As a consequence, integrating \eqref{2-16} from $1$ to $N$ we get
\[
c_N(w,\nu) =  -\frac{1}{\nu \Gamma(w+1)} + \int_1^N \frac{\xi^{\nu-1}}{\Gamma(w+1)} \Big(\big(1-\frac{1}{\xi}\big)^w-1\Big) \,\d \xi.
\]
Therefore, from \eqref{2-15} and \eqref{2-12} we obtain that for $N\geq 1$ and $u>0$  one has
\begin{equation}
\label{2-17}
T_N^{(1,3)}(w) = \frac{1}{\Gamma(w+1)} {\sum_{\substack{\nu=2 \\ 2\mid \nu}}^\infty}
 \Big\{ \int_1^N \big(\frac{\xi}{N}\big)^\nu \Big(\big(1-\frac{1}{\xi}\big)^w-1\Big)\, \frac{\d \xi}{\xi} 
- \frac{N^{-\nu}}{\nu}  - \frac{1}{\nu}  \Big(\big(1-\frac{1}{N}\big)^w-1\Big) \Big\}.
\end{equation}
But for $N\geq 2$ we have
\[ {\sum_{\substack{\nu=2 \\ 2\mid \nu}}^\infty}
 \frac{N^{-\nu}}{\nu} = \frac{1}{2} \sum_{\nu=1}^\infty \frac{N^{-2\nu}}{\nu} = -\frac{1}{2} \log \big(1-\frac{1}{N^2}\big)
\]
and
\[
\int_1^N\big(\frac{\xi}{N}\big)^\nu\, \frac{\d \xi}{\xi} + \frac{N^{-\nu}}{\nu} = \frac{1}{\nu},
\]
so that \eqref{2-17} becomes {
\[
\begin{split}
&= \frac{1}{\Gamma(w+1)} \Big\{ \frac{1}{2} \big(1-\frac{1}{N}\big)^w   \log \big(1-\frac{1} {N^2}\big)
 + 
 {\sum_{\substack{\nu=2 \\ 2\mid \nu}}^\infty}
 \int_1^N \big(\frac{\xi}{N}\big)^\nu \Big(\big(1-\frac{1}{\xi}\big)^w-\big(1-\frac{1}{N}\big)^w \Big)\, \frac{\d \xi}{\xi} \Big\} \\
&= \frac{1}{\Gamma(w+1)} \Big\{ \frac{1}{2} \big(1-\frac{1}{N}\big)^w \log \big(1-\frac{1} {N^2}\big)  + \int_1^N \Big(\big(1-\frac{1}{\xi}\big)^w-\big(1-\frac{1}{N}\big)^w \Big) \frac{\xi^2/N^2}{1-\xi^2/N^2}\,
\frac{\d\xi}{\xi} \Big\}.
\end{split}
\]}
Note that the above integral is convergent at $\xi=N$ since for $N\geq 2$ and $w$ in a compact {part } of $\CC$ we have
\begin{equation}
\label{2-18}
\big(1-\frac{1}{\xi}\big)^w-\big(1-\frac{1}{N}\big)^w \ll N-\xi
\end{equation}
{ uniformly in $w$ } as $\xi\to N^-$. Therefore, for $N\geq 2$ and $u>0$ we deduce that
{
\begin{equation}
\label{2-19}
T_N^{(1,3)}(w)= \frac1{\Gamma(w+1)} \Big\{  \frac12 \big(1-\frac{1}{N}\big)^w \log \big(1-\frac{1} {N^2}\big) 
 +   \int_1^N \Big(\big(1-\frac{1}{\xi}\big)^w-\big(1-\frac{1}{N}\big)^w \Big) \frac{\xi\,\d\xi}{N^2-\xi^2}\Big\}.
\end{equation}
}
Moreover, thanks to \eqref{2-18}, if $N\geq 3$ the part of the integral in \eqref{2-19} over $[2,N]$ extends to an entire function, and clearly the function
\[
 \int_1^2 \big(1-\frac{1}{N}\big)^w \frac{\xi\,\d\xi}{N^2-\xi^2}
\]
is also entire.

\smallskip
Consider now the remaining part of the second term on the right hand side of \eqref{2-19}, namely
\begin{equation}\label{neu}
I(w) = \frac{1}{\Gamma(w+1)}  \int_1^2 \big(1-\frac{1}{\xi}\big)^w \frac{\xi\,\d\xi}{N^2-\xi^2}.
\end{equation}
By the substitution $1-1/\xi=x$, and hence $\xi=1/(1-x)$ and $\d\xi=\d x/(1-x)^2$, we {find that}
\begin{equation}
\label{2-20}
I(w) =  \frac{1}{\Gamma(w+1)}  \int_0^{1/2} x^w \frac{\d x}{(1-x)(N^2(1-x)^2-1)}.
\end{equation}
Now we assume $N\geq 4$ and consider the function
\[
h_N(z) = \frac{1}{(1-z)(N^2(1-z)^2-1)} = \sum_{m=0}^\infty a_N(m) z^m,
\]
which is holomorphic in $|z|<1-1/N$. By Cauchy's {coefficient formula we see that whenever } $0<\delta<1-1/N$ we have $|a_N(m)| \leq \delta^{-m} \max_{|z|=\delta} |h_N(z)|$. {We take $\delta=2/3$ and infer that}
\begin{equation}
\label{2-21}
a_N(m) \leq 10 (3/2)^m.
\end{equation}
From \eqref{2-20} and \eqref{2-21} we therefore obtain that
\begin{equation}
\label{2-22}
I(w) = \frac{1}{\Gamma(w+1)}  \sum_{m=0}^\infty a_N(m) \int_0^{1/2} x^{w+m}\d x =  \frac{2^{-(w+1)}}{\Gamma(w+1)}  \sum_{m=0}^\infty \frac{a_N(m) 2^{-m}}{w+1+m}.
\end{equation}
By \eqref{2-21}, the last series converges absolutely,  and uniformly on {any  compact part of $\CC$ not containing any } of the points $w+1=-m$ with $m\geq0$, and hence it represents a meromorphic function with at most simple poles at $w+1=-m$, $m\geq0$. But such poles {cancel with } the zeros of $1/\Gamma(w+1)$, {so that } $I(w)$ is an entire function. Gathering the previous results of this subsection, we finally conclude that $T_N^{(1,3)}(w)$ is an entire function, for each $N\geq 4$.

\medskip
{\bf 2.6. Computing $T_N^{(1,4)}(w)$.} Since the term $T_N^{(1,4)}(w)$ is similar {to $T_N^{(1,3)}(w)$, a treatment along the lines of the previous subsection is possible. The details turn out to be } somewhat simpler, and hence we shall be more sketchy here. In analogy with \eqref{2-12}, \eqref{2-13} and \eqref{2-14}, for $N\geq 1$ and $u>0$ we have
\begin{equation}
\label{2-23}
T_N^{(1,4)}(w)  =  {\sum_{\substack{\nu=2 \\ 2\mid \nu}}^\infty}
 \left\{ \widetilde{H}_N(w,\nu)   + \frac{1}{\nu \Gamma(w+1)}  \Big(\big(1-\frac{1}{N}\big)^w-1\Big) \right\},
\end{equation}
where $ \widetilde{H}_N(w,\nu)$ satisfies
\begin{equation}
\label{2-24}
N^{-\nu}  \widetilde{H}_N(w,\nu) = \frac{1}{2\pi i} \int_{(-1/2)} \frac{N^{s-\nu}}{s-\nu}\, \frac{\Gamma(s)}{\Gamma(s+w+1)}\, \d s
\end{equation}
and
\[
\frac{\partial}{\partial N}  \widetilde{H}_N(w,\nu) - \frac{\nu}{N}  \widetilde{H}_N(w,\nu) = \frac{1}{N\Gamma(w+1)} \Big(\big(1-\frac{1}{N}\big)^w-1\Big). 
\]
Moreover, as in \eqref{2-15} we search for a function $\widetilde{c}_N(w,\nu)$ such that
\begin{equation}
\label{2-25}
\widetilde{H}_N(w,\nu) = \widetilde{c}_N(w,\nu) N^{\nu},
\end{equation}
which in analogy with \eqref{2-16} satisfies
\[
\frac{\partial}{\partial N} \widetilde{c}_N(w,\nu) = \frac{N^{-\nu-1}}{\Gamma(w+1)} \Big(\big(1-\frac{1}{N}\big)^w-1\Big).
\]
Now, from \eqref{2-24} and \eqref{2-25} we have that for any given $\nu\geq1$ and $w$ with $u>0$, as $N\to\infty$
\[
\widetilde{c}_N(w,\nu) \ll N^{-1/2-\nu},
\]
hence
\[
\widetilde{c}_N(w,\nu) = -\frac{1}{\Gamma(w+1)} \int_N^\infty  \Big(\big(1-\frac{1}{\xi}\big)^w-1\Big) \xi^{-\nu-1}\, \d \xi
\]
and therefore
\[
 \widetilde{H}_N(w,\nu)  =  -\frac{1}{\Gamma(w+1)} \int_N^\infty  \Big(\big(1-\frac{1}{\xi}\big)^w-1\Big) \big(\frac{N}{\xi}\big)^\nu\, \frac{\d \xi}{\xi}.
\]
Inserting this into \eqref{2-23}, in analogy with \eqref{2-19} we obtain that for $N\geq 2$ and $u>0$
\[
\begin{split}
T_N^{(1,4)}(w) &= \frac{1}{\Gamma(w+1)}  {\sum_{\substack{\nu=2 \\ 2\mid \nu}}^\infty}
  \Big\{ \Big(\big(1-\frac{1}{N}\big)^w-1\Big) \int_N^\infty \big(\frac{N}{\xi}\big)^\nu\, \frac{\d\xi}{\xi}
 - \int_N^\infty  \Big(\big(1-\frac{1}{\xi}\big)^w-1\Big) \big(\frac{N}{\xi}\big)^\nu\, \frac{\d \xi}{\xi} \Big\} \\
&= \frac{1}{\Gamma(w+1)} \int_N^\infty  \Big(\big(1-\frac{1}{N}\big)^w - \big(1-\frac{1}{\xi}\big)^w \Big) \frac{N/\xi}{1-N^2/\xi^2} \,\frac{\d\xi}{\xi} \\
&= \frac{N}{\Gamma(w+1)} \int_N^\infty  \Big(\big(1-\frac{1}{N}\big)^w - \big(1-\frac{1}{\xi}\big)^w \Big)\, \frac{\d\xi}{\xi^2-N^2}.
\end{split}
\]
Finally, since the last integral is uniformly convergent at $N$ {as long as $w$ ranges over a  compact part } of $\CC$, thanks to \eqref{2-18}, and clearly the same holds at $\infty$, the function
\begin{equation}
\label{2-26}
T_N^{(1,4)}(w) =  \frac{N}{\Gamma(w+1)} \int_N^\infty  \Big(\big(1-\frac{1}{N}\big)^w - \big(1-\frac{1}{\xi}\big)^w \Big)\, \frac{\d\xi}{\xi^2-N^2}
\end{equation}
is entire.

\medskip
The expression \eqref{1-9} in Proposition 1 and the fact that $T_N(w)$ is an entire function now follow gathering \eqref{2-3}, \eqref{2-5}, \eqref{2-7}, \eqref{2-10}, \eqref{2-11}, \eqref{2-19}, the conclusion of Subsection 2.5 and \eqref{2-26}.

\medskip
{\bf 2.7. Estimating $T_N(w)$.} In order to prove \eqref{1-10} we estimate the terms inside the brackets on the right hand side of  \eqref{1-9}. The first term is bounded by
\[
\ll \sum_{n\leq |w|} \frac{\Lambda(n)}{n} \Big(\big(1-\frac{1}{nN}\big)^u + 1\Big) + \sum_{n>|w|} \frac{\Lambda(n)}{n} \Big|\big(1-\frac{1}{nN}\big)^w -1\Big|,
\]
and, uniformly for $N\geq 4$, we have that
\[
\big(1-\frac{1}{nN}\big)^u + 1 \ll 2^{|u|} , \qquad  \big(1-\frac{1}{nN}\big)^w -1 = e^{w\log(1-1/(nN))}-1 \ll \frac{|w|}{n} \ \ \text{for} \ \ n>|w|.
\]
Therefore, by standard elementary bounds the above two sums are bounded by
\[
\ll 2^{|u|} \log(|w| + 2) \qquad \text{and} \qquad \ll |w| \sum_{n>|w|} \frac{\log n}{n^2} \ll \log(|w| + 2)
\]
respectively, and hence for any $w\in\CC$ and $N\geq 4$ one has
\begin{equation}
\label{2-27}
\sum_{n=1}^\infty \frac{\Lambda(n)}{n} \Big(\big(1-\frac{1}{nN}\big)^w-1\Big) \ll 2^{|u|} \log(|w|+2).
\end{equation}

\smallskip
The second term in \eqref{1-9} vanishes for $w=0$, and $(1-\xi/N)^w-1 \ll |w|\xi$ {holds } for $\xi<1/|w|$, for every $N\geq 4$. Hence such a term is $\ll 1$ when $|w|\leq 1$, while we split the range of integration into $[0,1/|w|] \cup [1/|w|,1]$ when $|w|>1$. In this case we have
\begin{equation}
\label{2-28}
\int_0^1\Big(\big(1-\frac{\xi}{N}\big)^w-1\Big)\, \frac{\d\xi}{\xi} \ll 1 +2^{|u|} \int_{1/|w|}^1 \frac{\d\xi}{\xi} \ll 2^{|u|} \log(|w|+2),
\end{equation}
and clearly the last bound in \eqref{2-28} holds for every $w\in\CC$ and $N\geq 4$.

\smallskip
Obviously, the third and fourth term in \eqref{1-9} are bounded, for every $w\in\CC$ and $N\geq4$, by
\begin{equation}
\label{2-29}
\ll 2^{|u|}.
\end{equation}

\smallskip
More care is required {for } the fifth term in \eqref{1-9}, {stemming from the integral in  \eqref{2-19}. In this integral we denote by $A(w)$  the part over $[1,2]$}. For $w\in\CC$ and $N\geq 4$ we have
\[
A(w) = \int_1^2 \big(1-\frac{1}{\xi}\big)^w\, \frac{\xi\,\d\xi}{N^2-\xi^2} + O\Big(2^{|u|} \int_1^2 \frac{\xi\,\d\xi}{N^2-\xi^2} \Big) =  A_1(w) + O\big(2^{|u|}\big),
\]
say. In view of {\eqref{neu} } and \eqref{2-22} we have
\[
A_1(w) = 2^{-(w+1)} \sum_{m=0}^\infty \frac{a_N(m) 2^{-m}}{w+1+m},
\]
and thanks to \eqref{2-21}, for $|w+m+1|>\delta$ we get $A_1(w) \ll 2^{|u|}/\delta$. Hence
\[
A(w) \ll \frac{2^{|u|}}{\delta}.
\]
In order to deal with the part over $[2,N]$ of the integral in \eqref{2-19} we note that
\[
\big(1-\frac{1}{\xi}\big)\big(1-\frac{1}{N}\big)^{-1} = 1 + O\big(\frac{|N-\xi|}{N}\big),
\]
and hence for $\xi> N(1-1/|w|)$ we have
\begin{equation}
\label{2-30}
\big(1-\frac{1}{\xi}\big)^w - \big(1-\frac{1}{N}\big)^w = \big(1-\frac{1}{N}\big)^w\Big(e^{w\log(1-1/\xi)(1-1/N)^{-1}} - 1\Big) \ll 2^{|u|} |w| \frac{|N-\xi|}{N}.
\end{equation}
{We temporarily assume that $|w|>(1-2/N)^{-1}$, and } split the integral over $[2,{N})$ into the part over $[2,N(1-1/|w|)]$ {and } the part over $[N(1-1/|w|),N]$, and denote by $B(w)$ and $C(w)$ {these } parts, respectively. A direct estimate gives
\[
B(w) \ll 2^{|u|} \int_2^{N(1-1/|w|)} \frac{\d\xi}{N-\xi} \ll 2^{|u|} \log(|w|+2),
\]
while thanks to \eqref{2-30} we obtain
\begin{equation}
\label{2-31}
C(w) \ll 2^{|u|} |w| \int_{N(1-1/|w|)}^N \frac{N-\xi}{N^2-\xi^2}\, \d\xi \ll 2^{|u|}.
\end{equation}
{The case where $|w|\leq(1-2/N)^{-1}$ is simpler. Here the whole integral over $[2,N]$ can be estimated as we bounded } $C(w)$. Therefore, gathering the above bounds for $A(w),B(w)$ and $C(w)$ we conclude that the  fifth term in \eqref{1-9} is bounded by
\begin{equation}
\label{2-32}
\ll \frac{2^{|u|} \log(|w|+2)}{\delta}
\end{equation}
for every $N\geq 4$ and $w\in\CC$ satisfying $|w+n+2|>\delta$.

\smallskip
Similarly, we split the integral over $[N,\infty)$ in the sixth term of \eqref{1-9} into the part over $[N,N(1+1/|w|)]$ plus the part over $[N(1+1/|w|), \infty)$, which we denote by $D(w)$ and $E(w)$, respectively. Arguing as in \eqref{2-31} we have
\[
N D(w) \ll 2^{|u|} |w| \int_N^{N(1+1/|w|)} \frac{\xi-N}{\xi^2-N^2}\, \d\xi \ll 2^{|u|},
\]
while a direct estimate gives
\[
NE(w) \ll 2^{|u|} N \int _{N(1+1/|w|)}^\infty \frac{\d\xi}{\xi^2-N^2}.
\]
Substituting $\xi-N =y$, decomposing $y^{-1}(y+2N)^{-1}$ into partial fractions and computing the resulting integrals we obtain
\[
\int _{N(1+1/|w|)}^\infty \frac{\d\xi}{\xi^2-N^2} = \frac{1}{2N} \log(2|w|+1).
\]
Hence the sixth term in \eqref{1-9} is bounded by
\begin{equation}
\label{2-33}
\ll 2^{|u|} \log(|w|+2)
\end{equation}
for every $N\geq 4$ and $w\in\CC$.

\medskip
Gathering \eqref{2-27},\eqref{2-28},\eqref{2-29},\eqref{2-32} and \eqref{2-33}, we obtain \eqref{1-10}, and the proof of Proposition 1 is complete. \fine

\bigskip
{\bf 2.8. Proof of Proposition 2.} Let $u>0$. {An application of \eqref{mellin}, followed by a shift of the line of integration to $\sigma=-1/2$, yields the identities}
\begin{equation}
\label{2-34}
\begin{split}
\frac{1}{\Gamma(w+1)} \sum_{n< N} \Lambda(n) \big(1-\frac{n}{N}\big)^w &= -\frac{1}{2\pi i} \int_{(2)} \frac{\zeta'}{\zeta}(s) \frac{\Gamma(s)}{\Gamma(s+w+1)} N^s\, \d s \\
&= \frac{N}{\Gamma(w+2)} - Z_N(w) - \frac{\zeta'}{\zeta}(0) \frac{1}{\Gamma(w+1)} + T_N(w).
\end{split}
\end{equation}
In view of Proposition 1, \eqref{2-34} gives the analytic continuation of $Z_N(w)$ to $\CC$. Moreover, for $|s+m| > \delta$ we have
\[
Z_N(w) \ll \frac{1}{|\Gamma(w+1)|}  \sum_{n<N} \Lambda(n) \big(1-\frac{n}{N}\big)^u + \frac{N}{\delta|\Gamma(w+1)|} + \frac{1}{|\Gamma(w+1)|} + |T_N(w)|.
\]
But clearly we have
\[
 \sum_{n<N} \Lambda(n) \big(1-\frac{n}{N}\big)^u \ll N^{|u|+1}
\]
for every $u\in\RR$, and for $u\leq -3/2$
\[
 \sum_{n<N} \Lambda(n) \big(1-\frac{n}{N}\big)^u \ll  N^{-u} \log N \sum_{n<N} n^u \ll N^{|u|} \log N.
\]
Proposition 2 now follows from the above bounds, thanks to \eqref{1-10}. \fine

\bigskip
\section{Proof of the theorem}

\smallskip

{\bf 3.1. The first shift.} Our point of departure is \eqref{3-2} where we 
 shift the $s$-integration  to the line $\sigma=-1/2$. We recall \eqref{1-7} and \eqref{1-8}, thus getting
\begin{equation}
\label{3-3}
\begin{split}
G_k(N) &= \frac{1}{2\pi i} \int_{(2)} -\frac{\zeta'}{\zeta}(w) \frac{\Gamma(w)}{\Gamma(w+k+2)} N^{w+1}\, \d w  -  \frac{1}{2\pi i} \int_{(2)} -\frac{\zeta'}{\zeta}(w) \Gamma(w) Z_N(w+k) N^w\, \d w \\
&\hskip1cm - \frac{\zeta'}{\zeta}(0) \frac{1}{2\pi i} \int_{(2)} -\frac{\zeta'}{\zeta}(w) \frac{\Gamma(w)}{\Gamma(w+k+1)} N^w\, \d w \\\
&\hskip1cm + \frac{1}{2\pi i} \int_{(2)} -\frac{\zeta'}{\zeta}(w) \Gamma(w) T_N(w+k) N^w\, \d w \\
&= N\Gamma(N,k+1) - Z(N,k) - \frac{\zeta'}{\zeta}(0) \Gamma(N,k) + T(N,k),
\end{split}
\end{equation}
say. All integrals here are absolutely convergent for $k>0$, thanks to Stirling's formula and the bounds \eqref{1-10} and \eqref{1-11}. 

\medskip
{\bf 3.2. Shifting to $-\infty$.} Next, we want to shift to $-\infty$ {all } $w$-integrations in \eqref{3-3}. To this end, we first need to get suitable bounds for the integrands. We may clearly assume that $|w|$ is sufficiently large and $|w+m|\geq 1/4$ for every integer $m\geq1$. Recalling \eqref{2-1}, \eqref{2-2} and the reflection formula for the $\Gamma$ function, we have
\begin{equation}
\label{3-4}
\frac{\zeta'}{\zeta}(w) \ll \left|\frac{g'}{g}(w) \right|+ 1 \ll \left| \frac{\Gamma'}{\Gamma}(1-w)  \right| + \left| \frac{\sin(\pi w/2)}{\cos(\pi w/2)} \right| + \left|  \frac{\cos(\pi w)}{\sin(\pi w)} \right| +1 \ll \log|w|,
\end{equation}
thanks to (6) in Chapter 10 of \cite{Dav/1980} and the bounds $O(1)$ for the two trigonometric terms. In view of the shape of the bounds \eqref{1-10} and \eqref{1-11}, we apply again the reflection formula and then Stirling's formula to bound $\Gamma(w)/\Gamma(w+k+1)$ {to conclude } that
\begin{equation}
\label{3-5}
\frac{\Gamma(w)}{\Gamma(w+k+1)} \ll |w|^{-(k+1)}.
\end{equation}

\smallskip
Since the computations involved in the shift to $-\infty$ of the integrals in \eqref{3-3} are now quite standard, we only give a sketch of the argument. We treat explicitly only the integral $Z(N,k)$, since $T(N,k)$ is similar but easier and $\Gamma(N,k)$ gives rise to a classical weighted explicit formula. We first restrict the integration on the line $u=-1/2$ to the segment with $|v|\leq V$ and then shift such a segment to $u=-(U+1/2)$, where $0<U<V$ are sufficiently large and $U\in\NN$. We denote by $Z_{U,V}^{\text{hor}}(N,k)$ the integral over the two horizontal sides $[-(U+1/2)\pm iV,-1/2 \pm iV]$ and by $Z_{U,V}^{\text{vert}}(N,k)$ the one over the vertical side $[-(U+1/2)-iV,-(U+1/2)+iV]$. Thanks to \eqref{1-11} with $\delta=1/4$, \eqref{3-4} and \eqref{3-5} we have
\[
Z_{U,V}^{\text{hor}}(N,k) \ll N^cU \frac{\log^2V}{V^{k+1}} \quad \text{and} \quad Z_{U,V}^{\text{vert}}(N,k) \ll N^c \frac{\log^2U}{U^k}
\]
for some $c>0$, and hence
\[
\lim_{U\to+\infty} \big(\lim_{V\to+\infty} Z_{U,V}^{\text{hor}}(N,k) + Z_{U,V}^{\text{vert}}(N,k) \big) =0.
\]
Therefore, $Z(N,k)$ equals the sum of the residues { $-\frac{\zeta'}{\zeta}(w) \Gamma(w) Z_N(w+k)N^w$}. Clearly, {this argument also applies to } the other three integrals in \eqref{3-3}.

\medskip
{\bf 3.3. Computing the residues.} For each of the four integrals in \eqref{3-3} we have to compute the following residues:

\noindent
(a) at the simple pole of $-\zeta'/\zeta(w)$ at $w=1$;

\noindent
(b) at the simple poles of $-\zeta'/\zeta(w)$ at $w=\rho$, $\rho$ non-trivial zero;

\noindent
(c) at the simple poles of $\Gamma(w)$ at $w=0$ and $w=-\nu$, $\nu\geq 1$ and $2\nmid \nu$;

\noindent
(d) at the double poles of $-\frac{\zeta'}{\zeta}(w) \Gamma(w)$ at $w=-\nu$, $\nu\geq 1$ and $2|\nu$.

\smallskip
The residues of type (a) produce the terms
\begin{equation}
\label{3-6}
\frac{N^2}{\Gamma(k+3)} - N Z_N(k+1) -\frac{\zeta'}{\zeta}(0) \frac{N}{\Gamma(k+2)} +N T_N(k+1),
\end{equation}
while the type (b) residues give rise to the sums
\begin{equation}
\label{3-7}
\begin{split}
&-\sum_\rho \frac{\Gamma(\rho)}{\Gamma(\rho+k+2)} N^{\rho+1} + \sum_\rho \Gamma(\rho) Z_N(\rho+k) N^\rho \\
&+\frac{\zeta'}{\zeta}(0) \sum_\rho \frac{\Gamma(\rho)}{\Gamma(\rho+k+1)} N^\rho -  \sum_\rho \Gamma(\rho) T_N(\rho+k) N^\rho.
\end{split}
\end{equation}
By Stirling's formula, the sums in \eqref{3-7} are absolutely convergent thanks to \eqref{1-10} and \eqref{1-11}, since $k>0$. The residues of type (c) produce the terms
\begin{equation}
\label{3-8}
-\frac{\zeta'}{\zeta}(0) \frac{N}{\Gamma(k+2)} + \frac{\zeta'}{\zeta}(0) Z_N(k) + \frac{\zeta'}{\zeta}(0)^2 \frac{1}{\Gamma(k+1)} -\frac{\zeta'}{\zeta}(0) T_N(k)
\end{equation}
plus the sums
\begin{equation}
\label{3-9}
\begin{split}
&\sum_{\substack{\nu\geq 1 \\ 2\nmid \nu}} \frac{\zeta'}{\zeta}(-\nu) \frac{N^{-\nu+1}}{\nu! \Gamma(-\nu+k+2)} -  \sum_{\substack{\nu\geq 1 \\ 2\nmid \nu}} \frac{\zeta'}{\zeta}(-\nu) Z_N(-\nu+k) \frac{N^{-\nu}}{\nu!} \\
&-\frac{\zeta'}{\zeta}(0) \sum_{\substack{\nu\geq 1 \\ 2\nmid \nu}} \frac{\zeta'}{\zeta}(-\nu) \frac{N^{-\nu}}{\nu! \Gamma(-\nu+k+1)} + \sum_{\substack{\nu\geq 1 \\ 2\nmid \nu}} \frac{\zeta'}{\zeta}(-\nu) T_N(-\nu+k) \frac{N^{-\nu}}{\nu!}.
\end{split}
\end{equation}
Again, the sums in \eqref{3-9} are absolutely convergent for $k>0$, thanks to \eqref{1-10}, \eqref{1-11}, \eqref{3-4} and Stirling's formula.

\smallskip
A bit more care is required {for } residues of type (d). For even integers $\nu\geq 1$, let $a_\nu$ and $b_\nu$ be defined by the following Laurent expansions at $w=-\nu$ {via}
\begin{equation}
\label{3-10}
-\frac{\zeta'}{\zeta}(w) = \frac{1}{w+\nu} +a_\nu + \dots \qquad  \Gamma(w) = \frac{1}{\nu!(w+\nu)} + \frac{b_\nu}{\nu!} +\dots,
\end{equation}
and let
\begin{equation}
\label{3-11}
A_\nu(N) = 
\begin{cases}
a_\nu+b_\nu+\log N & \text{if $2|\nu$} \\
-\frac{\zeta'}{\zeta}(-\nu) & \text{if $2\nmid \nu$.}
\end{cases}
\end{equation}
Moreover, denoting by $F(w)$ any of the four functions 
\[
\frac{N^{w+1}}{\Gamma(w+k+2)},\quad Z_N(w+k)N^w,\quad \frac{N^w}{\Gamma(w+k+1)}, \quad T_N(w+k)N^w,
\] 
in view of \eqref{3-10} the residues of type (d) at $w=-\nu$ are of the form
\[
\frac{a_\nu+b_\nu}{\nu!} F(-\nu) + \frac{1}{\nu!} F'(-\nu).
\]
Hence by \eqref{3-11} the contribution of the residues of type (d) is
\begin{equation}
\label{3-12}
\begin{split}
&\sum_{\substack{\nu\geq 1 \\ 2| \nu}} \frac{A_\nu(N)}{\nu!} \frac{N^{-\nu+1}}{\Gamma(-\nu+k+2)} - \sum_{\substack{\nu\geq 1 \\ 2| \nu}}  \frac{\Gamma'(-\nu+k+2)}{\Gamma^2(-\nu+k+2)} \frac{N^{-\nu+1}}{\nu!} \\
& \hskip.5cm - \sum_{\substack{\nu\geq 1 \\ 2| \nu}} \frac{A_\nu(N)}{\nu!} Z_N(k-\nu) N^{-\nu} -  \sum_{\substack{\nu\geq 1 \\ 2| \nu}}  Z'_N(k-\nu) \frac{N^{-\nu}}{\nu!} \\
&\hskip.5cm -\frac{\zeta'}{\zeta}(0) \sum_{\substack{\nu\geq 1 \\ 2| \nu}} \frac{A_\nu(N)}{\nu!} \frac{N^{-\nu}}{\Gamma(-\nu+k+1)} + \frac{\zeta'}{\zeta}(0) \sum_{\substack{\nu\geq 1 \\ 2| \nu}} \frac{\Gamma'(-\nu+k+1)}{\Gamma^2(-\nu+k+1)} \frac{N^{-\nu}}{\nu!} \\
&\hskip.5cm+  \sum_{\substack{\nu\geq 1 \\ 2| \nu}} \frac{A_\nu(N)}{\nu!} T_N(k-\nu) N^{-\nu} +  \sum_{\substack{\nu\geq 1 \\ 2| \nu}} T'_N(k-\nu) \frac{N^{-\nu}}{\nu!}.
\end{split}
\end{equation}
The absolute convergence of the second and sixth sum in \eqref{3-12} follows by computations similar to those leading to \eqref{3-4} and \eqref{3-5}, giving in particular that for $a>0$
\begin{equation}
\label{3-13}
\frac{1}{\nu! \Gamma(-\nu + a)} \ll \frac{1}{\nu^a}.
\end{equation}
Hence now we concentrate on the remaining six sums. However, it is already clear from \eqref{3-3}, \eqref{3-6}-\eqref{3-8} after a simple a rearrangement of terms, \eqref{3-9} and \eqref{3-12} that the explicit formula in the theorem holds as a formal identity. 

\medskip
{\bf 3.4. The final estimates.} We need the following lemma.

\medskip
{\bf Lemma 3.1.} {\sl For $N\geq4$ and $\nu\geq 1$ we have}
\[
A_\nu(N) \ll \log(\nu N).
\]

\medskip
{\it Proof.} The lemma follows at once from \eqref{3-4} and \eqref{3-11} if $2\nmid \nu$. If $2|\nu$, then by \eqref{3-10} we have
\[
a_\nu = -\frac{1}{2\pi i} \int_{|w+\nu|=1} \big(\frac{\zeta'}{\zeta}(w) + \frac{1}{w+\nu}\big) \d w \ll \log (\nu+1)
\]
thanks to \eqref{3-4}, and
\[
\begin{split}
b_\nu &= \nu! \lim_{w\to-\nu} \big(\Gamma(w) - \frac{1}{\nu!(w+\nu)}\big) = \nu!  \lim_{w\to-\nu} \Big(\frac{\Gamma(w)(w+\nu)-1/\nu!}{w+\nu}\Big) \\
&=  \nu! \lim_{w\to-\nu} \big(\Gamma(w) + (w+\nu)\Gamma'(w)\big)  = \nu! \lim_{w\to-\nu} \big((w+\nu)\Gamma(w)\big) \Big(\frac{1}{w+\nu} + \frac{\Gamma'}{\Gamma}(w)\Big) \\
& =  \lim_{w\to-\nu} \Big( \frac{\Gamma'}{\Gamma}(w) + \frac{1}{w+\nu}\Big),
\end{split}
\]
since $2|\nu$. Taking the logarithmic derivative of the Hadamard product of $1/\Gamma(w)$ as in \eqref{2-6} we see that
\[
\frac{\Gamma'}{\Gamma}(w) = -\gamma -\frac{1}{w} + \sum_{n=1}^\infty \Big(\frac{1}{n} - \frac{1}{w+n}\Big),
\]
hence
\[
b_\nu = -\gamma + \frac{2}{\nu} + \sum_{n\neq \nu} \Big(\frac{1}{n} - \frac{1}{n-\nu}\Big) \ll \log(\nu+1).
\]
The lemma now follows thanks to \eqref{3-11}. \fine

\medskip
Now we prove the absolute convergence of the above mentioned six sums in \eqref{3-12}. Thanks to Lemma 3.1, the absolute convergence of the fifth sum in \eqref{3-12} follows from \eqref{3-13}. Next we deal with the two sums involving the function $Z_N(w)$, those involving $T_N(w)$ being similar and easier; we may clearly assume that $\nu$ is sufficiently large. If $k$ is not an integer we may apply directly the second bound for $Z_N(k-\nu)$ in \eqref{1-11}, choosing any $0<\delta<k$. If $k$ is an integer, by the maximum principle we have that
\[
|Z_N(k-\nu)| \leq \max_{|w+k-\nu|=\delta} |Z_N(w)|
\]
with a small $0<\delta<k$, and then we apply again the second bound in \eqref{1-11}. Hence in both cases we have
\begin{equation}
\label{3-14}
Z_N(k-\nu) \ll \frac{N^{\nu-k+\delta} \log N + 2^{\nu-k+\delta} \log \nu}{|\Gamma(-\nu+k-\delta+1)|}.
\end{equation}
The bound for $Z'_N(k-\nu)$ is obtained in a similar way, using of Cauchy's formula
\[
Z'_N(k-\nu) = \frac{1}{2\pi i} \int_{|w+\nu-k|=\delta} \frac{Z_N(w)}{(w+\nu-k)^2} \d w
\]
and hence getting that
\begin{equation}
\label{3-15}
Z'_N(k-\nu) \ll \frac{N^{\nu-k+\delta} \log N + 2^{\nu-k+\delta} \log \nu}{|\Gamma(-\nu+k-\delta+1)|}.
\end{equation}
with any $0<\delta<k$. The absolute convergence of the third and fourth sum in \eqref{3-12} follows now from Lemma 3.1, \eqref{3-13}, \eqref{3-14} and \eqref{3-15}, since $0<\delta<k$. The theorem is therefore proved.

\medskip
Finally, from \eqref{1-12}-\eqref{1-14}, \eqref{3-9} and \eqref{3-12} we have that
\begin{equation}
\label{3-16}
\Sigma_\Gamma(N,k) = \sum_{\substack{\nu\geq 1 \\ 2\nmid \nu}}  \frac{\zeta'}{\zeta}(-\nu) \frac{N^{-\nu}}{\nu! \Gamma(-\nu+k+1)} + \sum_{\substack{\nu\geq 1 \\ 2| \nu}} \Big(\frac{A_\nu(N)}{\Gamma(-\nu+k+1)} -  \frac{\Gamma'(-\nu+k+2)}{\Gamma^2(-\nu+k+1)} \Big) \frac{N^{-\nu}}{\nu!},
\end{equation}
\begin{equation}
\label{3-17}
\Sigma_Z(N,k) = \sum_{\substack{\nu\geq 1 \\ 2\nmid \nu}} \frac{\zeta'}{\zeta}(-\nu) Z_N(-\nu+k) \frac{N^{-\nu}}{\nu!} +\sum_{\substack{\nu\geq 1 \\ 2| \nu}} \Big(A_\nu(N) Z_N(k-\nu) +  Z'_N(k-\nu) \Big) \frac{N^{-\nu}}{\nu!},
\end{equation}
\begin{equation}
\label{3-18}
\Sigma_T(N,k) = \sum_{\substack{\nu\geq 1 \\ 2\nmid \nu}} \frac{\zeta'}{\zeta}(-\nu) T_N(-\nu+k) \frac{N^{-\nu}}{\nu!} +\sum_{\substack{\nu\geq 1 \\ 2| \nu}} \Big(A_\nu(N) T_N(k-\nu) +  T'_N(k-\nu) \Big) \frac{N^{-\nu}}{\nu!};
\end{equation}
moreover, we denote by $\Sigma_\Gamma^M(N,k)$, $\Sigma_Z^M(N,k)$ and $\Sigma_T^M(N,k)$ the tail of such sums, from $\nu=M+1$ to $\infty$. Hence from \eqref{3-4}, \eqref{3-13}, Lemma 3.1, \eqref{3-14} and \eqref{3-15}, and the analogues of  \eqref{3-14} and \eqref{3-15} for $T_N(k-\nu)$ where the term $N^{\nu-k+\delta}\log N$ is missing, {it follows } that cutting the above sums at $\nu=M\geq2$ produces the errors
\begin{equation}
\label{3-19}
\begin{split}
\Sigma_\Gamma^M(N,k) &\ll \frac{N^{-(M+1)}\log (NM)}{M^{k}}, \\
\Sigma_Z^M(N,k) &\ll_\delta  \frac{N^{-k+\delta} \log^2(NM)}{M^{k-\delta}},  \\
\Sigma_T^M(N,k) &\ll_\delta  \frac{(N/2)^{-(M+1)} \log^2(NM)}{M^{k-\delta}},
\end{split}
\end{equation}
for every $\delta>0$. All the statements are now proved.


\ifx\undefined\bysame{poly}.
\newcommand{\bysame}{\leavevmode\hbox to3em{\hrulefill}\ ,}
\fi

\vskip 1cm
\noindent
J\"org Br\"udern, Mathematisches Institut, Bunsenstr. 3-5, 37073 G\"ottingen, Germany. e-mail: bruedern@uni-math.gwdg.de

\medskip
\noindent
Jerzy Kaczorowski, Faculty of Mathematics and Computer Science, A.Mickiewicz University, 61-614 Pozna\'n, Poland and Institute of Mathematics of the Polish Academy of Sciences, 
00-956 Warsaw, Poland. e-mail: kjerzy@amu.edu.pl

\medskip
\noindent
Alberto Perelli, Dipartimento di Matematica, Universit\`a di Genova, via Dodecaneso 35, 16146 Genova, Italy. e-mail: perelli@dima.unige.it

\end{document}